%% file: main.tex
\documentclass[12pt]{amsart}

\usepackage[margin = 1in]{geometry}
\usepackage{amsfonts, amsmath,amscd, amssymb,
latexsym, mathrsfs, mathtools, 
slashed, stmaryrd, verbatim, wasysym, bbm, url}
\usepackage[usenames]{xcolor}
\usepackage{enumerate}
\usepackage[shortlabels]{enumitem}
\usepackage[utf8]{inputenc}

\usepackage[backend=bibtex,bibstyle=authoryear,citestyle=alphabetic,firstinits=true,doi=false,isbn=false,url=false,backref]{biblatex}
\DeclareFieldFormat{labelalphawidth}{\mkbibbrackets{#1}}

\defbibenvironment{bibliography}
  {\list
     {\printtext[labelalphawidth]{%
        \printfield{prefixnumber}%
    \printfield{labelalpha}%
        \printfield{extraalpha}}}
     {\setlength{\labelwidth}{\labelalphawidth}%
      \setlength{\leftmargin}{\labelwidth}%
      \setlength{\labelsep}{\biblabelsep}%
      \addtolength{\leftmargin}{\labelsep}%
      \setlength{\itemsep}{\bibitemsep}%
      \setlength{\parsep}{\bibparsep}}%
      }
  {\endlist}
  {\item}

\newtheorem*{acknowledgements}{Acknowledgements}
\newtheorem*{theorem*}{Theorem}
\newtheorem{theorem}{Theorem}
\newtheorem{corollary}{Corollary}

\theoremstyle{definition}
\newtheorem{example}{Example}
\newtheorem{remark}[example]{Remark}

\numberwithin{equation}{section}

\let\oldsqrt\sqrt
\def\sqrt{\mathpalette\DHLhksqrt}
\def\DHLhksqrt#1#2{%
\setbox0=\hbox{$#1\oldsqrt{#2\,}$}\dimen0=\ht0
\advance\dimen0-0.2\ht0
\setbox2=\hbox{\vrule height\ht0 depth -\dimen0}%
{\box0\lower0.4pt\box2}}

\usepackage{verbatim} 
\DeclareFontFamily{U}{mathx}{\hyphenchar\font45}
\DeclareFontShape{U}{mathx}{m}{n}{
      <5> <6> <7> <8> <9> <10>
      <10.95> <12> <14.4> <17.28> <20.74> <24.88>
      mathx10
      }{}
\DeclareSymbolFont{mathx}{U}{mathx}{m}{n}
\DeclareFontSubstitution{U}{mathx}{m}{n}
\DeclareMathAccent{\widecheck}{0}{mathx}{"71}

\renewcommand{\tilde}{\widetilde}

\renewcommand{\hat}[1]{\widehat{#1}}

\newcommand\eps\varepsilon
\renewcommand\epsilon\varepsilon

\newcommand{\abs}[1]{\left\lvert #1 \right\rvert}
\newcommand{\smallabs}[1]{\lvert #1 \rvert}

\newcommand\floor[1]{\lfloor #1 \rfloor}

\renewcommand\Im{\operatorname{Im}}

\newcommand\Mand{\text{ and }}

\newcommand\paperintro%
        {%
         }
\newcommand\paperbody%
        {%
         }

\newcommand\bbC{\mathbb{C}}

\newcommand\bbR{\mathbb{R}}

\newcommand\bbZ{\mathbb{Z}}

\newcommand\cF{\mathcal{F}}

\newcommand\cO{\mathcal{O}}

\DeclareMathAlphabet{\mathpzc}{OT1}{pzc}{m}{it}

\newcommand{\sbs}{\subset}

\usepackage[refpage]{nomencl}

\makenomenclature

\usepackage{tikz, tikz-cd, tkz-euclide}
\usetikzlibrary{calc,positioning,intersections,quotes,decorations.markings}
\makeatletter
\def\@tocline#1#2#3#4#5#6#7{\relax
  \ifnum #1>\c@tocdepth 
  \else
    \par \addpenalty\@secpenalty\addvspace{#2}%
    \begingroup \hyphenpenalty\@M
    \@ifempty{#4}{%
      \@tempdima\csname r@tocindent\number#1\endcsname\relax
    }{%
      \@tempdima#4\relax
    }%
    \parindent\z@ \leftskip#3\relax \advance\leftskip\@tempdima\relax
    \rightskip\@pnumwidth plus4em \parfillskip-\@pnumwidth
    #5\leavevmode\hskip-\@tempdima
      \ifcase #1
       \or\or \hskip 1em \or \hskip 2em \else \hskip 3em \fi%
      #6\nobreak\relax
    \hfill\hbox to\@pnumwidth{\@tocpagenum{#7}}\par
    \nobreak
    \endgroup
  \fi}
\makeatother

\def\annu#1{_{%
  \vbox{\hrule height .2pt 
    \kern 1pt 
    \hbox{$\scriptstyle {#1}\kern 1pt$}%
  }\kern-.05pt 
  \vrule width .2pt 
}}

\allowdisplaybreaks
\usepackage[colorlinks,citecolor=blue,urlcolor=red,bookmarks=false,hypertexnames=true]{hyperref} 
\hypersetup{
  colorlinks,
  citecolor=blue,
  linkcolor=blue,
  urlcolor=blue}
  
\title[Weyl asymptotics for functional difference operators]{Weyl asymptotics for functional difference operators with power to quadratic exponential potential}
\author{Yaozhong Qiu}
\address{Department of Mathematics, Imperial College London, 180 Queen’s
Gate, London, SW7 2AZ, United Kingdom}
\email{y.qiu20@imperial.ac.uk}

\begin{document}

\begin{abstract}
We continue the program first initiated in \cite{LST15} and develop a modification of the technique introduced in that paper to study the spectral asymptotics, namely the Riesz means and eigenvalue counting functions, of functional difference operators $H_0 = \cF^{-1} M_{\cosh(\xi)} \cF$ with potentials of the form $W(x) = \abs{x}^pe^{\abs{x}^\beta}$ for either $\beta = 0$ and $p > 0$ or $\beta \in (0, 2]$ and $p \geq 0$. We provide a new method for studying general potentials which includes the potentials studied in \cite{LST15, LST19}. The proof involves dilating the variance of the gaussian defining the coherent state transform in a controlled manner preserving the expected asymptotics. 
\end{abstract}

\maketitle

\input{introduction}

\input{main-results}

\begin{acknowledgements}
\textup{The author would like to thank Ari Laptev for introducing the author to the problem, his warm encouragement, and for helpful discussions. We are also very grateful to the anonymous referee for their thorough reading and comments which helped us improve the presentation of the paper. The author is supported by the President's Ph.D. Scholarship of Imperial College London.}
\end{acknowledgements}

\printbibliography

\end{document}

%% file: introduction.tex
\section{Introduction and main results}
In this paper, we continue the program first initiated in \cite{LST15} and study the Riesz means of some classes of functional difference operators with potential $H = H_0 + W$ on the line using the coherent state transform. In \cite{LST15} the functional difference operator $H_0$ is defined as the Fourier multiplier $2\cosh(2\pi b\xi)$ for $b > 0$, i.e.
\[ (H_0\psi, \varphi) = \int_\bbR (H_0\psi(x))\overline{\varphi(x)}dx = \int_\bbR 2\cosh(2\pi b\xi)\hat{\psi}(\xi)\overline{\hat{\varphi}(\xi)}d\xi = (2\cosh(2\pi b\xi)\hat{\psi}, \hat{\varphi}) \]
for all $\psi, \varphi \in L^2(\bbR)$ such that $2\cosh(2\pi b\xi)\hat{\psi} \in L^2(\bbR)$ with $\hat{\psi}(\xi) = \int_\bbR e^{-2\pi ix\xi}\psi(x)dx$ the Fourier transform of $\psi$. Equivalently, its action may be described as the sum of two complex shift operators
\[ (H_0\psi)(x) = \psi(x + ib) + \psi(x - ib) \]
for all $\psi \in L^2(\bbR)$ equipped with an analytic continuation to the strip $\{z \in \bbC \mid \abs{\Im z} < b\} \sbs \bbC$ such that ${\psi(- + iy) \in L^2(\bbR)}$ for all $\abs{y} < b$ and ${\psi(x \pm ib) \vcentcolon= \lim_{\epsilon \to 0^+} \psi(x \pm ib \mp i\epsilon)}$ defined in the sense of $L^2(\bbR)$-convergence. It was shown in \cite[Proposition~2.1]{LST15} that the unbounded operator $H$, for $W$ a continuous real-valued function blowing up at infinity in the sense $\inf_{\abs{x} \geq R} W(x) \rightarrow \infty$ as $R \rightarrow \infty$, admits a selfadjoint extension to $L^2(\bbR)$, which we continue to denote by $H$, with pure discrete spectrum. Denoting $(\lambda_j)_{j \geq 1} \sbs \bbR_{\geq 0}$ the nondecreasing sequence of eigenvalues of $H$, it was shown in \cite[Theorem~2.2]{LST15} the Riesz mean $\sum_{j \geq 1} (\lambda - \lambda_j)_+$ of $H = H_0 + 2\cosh(2\pi bx)$ has asymptotics
\[ \sum_{j \geq 1} (\lambda - \lambda_j)_+ = \frac{\lambda\log^2\lambda}{(\pi b)^2} + \cO(\lambda\log\lambda) \]
as $\lambda \rightarrow \infty$. Note the leading term coincides with the leading term of the classical phase space integral 
\[ \int_{\bbR^2} (\lambda - 2\cosh(2\pi bk) - 2\cosh(2\pi by))_+dkdy \]
as $\lambda \rightarrow \infty$ and where $2\cosh(2\pi bk) + 2\cosh(2\pi by)$ is the total symbol of $H_0 + 2\cosh(2\pi bx)$, and so the results of \cite{LST15} are Weyl type asymptotics (after the original result \cite{Wey11} by Weyl, see also \cite{Ivr16} for some history) linking the asymptotics of quantum mechanical expressions to classical phase space integrals. 

To our best knowledge, the literature on the study of the Riesz mean of functional difference operators using the coherent state transform is presently limited to the original paper \cite{LST15}, which studied the more general potential $W(x) = e^{2\pi bx} + \zeta e^{-2\pi bx}$ for $\zeta > 0$, whose asymptotics have the same leading term as the special case $\zeta = 1$ described above, and the companion paper \cite{LST19} which studied $W(x) = x^{2n} + r(x)$ for $n \in \bbZ_{\geq 1}$ and a continuous real-valued function $r$ satisfying $\abs{r(x)} \leq C(1 + \abs{x}^{2N-\epsilon})$ for $0 < \epsilon \leq 2N-1$. Functional difference operators have been studied in the literature prior to these papers without, however, using the coherent state transform. The case $W(x) = e^{2\pi bx}$, which has absolutely continuous spectrum $[2, \infty)$, was studied in \cite{Kas01a, Kas01b, TF15} and a generalisation thereof in \cite{Mic15}. Moreover, the coherent state transform has appeared elsewhere in the literature (that is not necessarily in the context of functional difference operators), for instance in \cite{Gei14} and in \cite{LW21} where it was used to study the fractional and logarithmic laplacians respectively, the former generalising \cite{Wey11}. A more general discussion of Riesz means and the coherent state transform can be found, for instance, in \cite{HH11} and \cite[Chapter~12]{LL01} respectively. 

For particular choices of $W$, the operator $H$ has physical meaning. The operator $H_0 + e^{2\pi bx}$ first appeared in \cite{TF86} and later in the representation theory of the quantum group $\text{SL}_q(2, \bbR)$, while the operator $H_0 + 2\cosh(2\pi bx)$ first appeared in \cite{Aga06} in the study of local mirror symmetry of toric Calabi-Yau manifolds. Denoting $U\psi(x) = \psi(x + ib)$ and $V = e^{2\pi bx}$, it was discovered in \cite{Aga06} functional difference operators can be regarded as quantisations of algebraic curves, namely both $H_0 + e^{2\pi bx}$ and $H_0 + 2\cosh(2\pi bx)$ can be realised as the quantisation of ${e^x + e^{-x} + e^y}$ and ${e^x + e^{-x} + e^y + e^{-y}}$ respectively. Spectral properties of these operators, together with the operator $H_{m, n} = U + V + q^{-mn}U^{-m}V^{-n}$, were also studied in \cite{GHM16, KM16, DKB21} and \cite[\S3]{LST15}. The operator $H_0 + x^{2n}$ and more generally with $x^{2n}$ replaced by any polynomial (and possibly of odd degree, and thus unbounded below) appeared in the recent work \cite{GM19} but also as early as \cite{PG96}. For more details on the relationship between spectral theory with mirror symmetry and string theory, we refer the reader to the survey \cite{marino2018spectral} and references therein. 

The principle behind the coherent state transform (in nice enough settings) is that the action of $H$ can be expressed as a multiplication operator in ``coherent state space" in analogue with the intuition that a pseudodifferential operator is a multiplication operator in phase space. Let $a > 0$ and $g_a(y) = (a/\pi)^{1/4}e^{-ay^2/2}$ be a $L^2(\bbR)$-normalised gaussian. We define the coherent state transform of $\psi \in L^2(\bbR)$ with respect to $g_a$ by the formula 
\begin{equation}\label{coherentstateformula}
\tilde{\psi}(k, y) \vcentcolon= \int_\bbR e^{-2\pi ikx}g_a(x-y)\psi(x)dx.
\end{equation}
It satisfies the formulae \cite[pp.~5]{LST15}
\begin{align}
\int_\bbR \smallabs{\tilde{\psi}(k, y)}^2dk &= (\abs{\psi}^2 \ast \abs{g_a}^2)(y), \label{eq6lst15} \\ 
\int_\bbR \smallabs{\tilde{\psi}(k, y)}^2dy &= (\smallabs{\hat{\psi}}^2 \ast \smallabs{\hat{g_a}}^2)(y), \label{eq7lst15}
\end{align}
which can be proved by the Plancherel and convolution theorems for the Fourier transform. More important to this paper is the formula \cite[\S12.12]{LL01}
\begin{equation}\label{lst19formula}
\int_{\bbR^2} W_0(y)\smallabs{\tilde{\psi}(k, y)}^2dkdy = \int_\bbR (W_0 \ast g_a^2)(x)\abs{\psi(x)}^2dx
\end{equation}
which expresses the fact multiplication by $(k, y) \mapsto W_0(y)$ in coherent state space $\bbR^2 = \bbR_k \times \bbR_y$ amounts to multiplication by ${W_0 \ast g_a^2}$ in spatial coordinates. 

The problem of finding the $W_0$ for which ${W_0 \ast g_a^2} = W$, that is finding the multiplication operator $W_0(y)$ in coherent state space corresponding to $W(x)$ in spatial coordinates, is generally difficult. By regarding convolution with $g_a^2$ as convolution with the euclidean heat kernel for some time $t > 0$ depending on $a$, it follows solving the convolution equation ${W_0 \ast g_a^2} = W$ for $W_0$ is essentially a matter of solving the backward heat equation, which is generally ill-posed. 

For some simple examples this can be done; for instance if $W(x) = \cosh(2\pi bx)$ then $W_0(y) = C_1\cosh(2\pi by)$ for some constant $C_1 > 0$ depending only on $a$ and $b$ \cite[pp.~6]{LST15}, and if $W(x) = x^{2n}$ for $n \in \bbZ_{\geq 1}$ then $W_0(y) = y^{2n} + R(y)$ where $R$ is a polynomial of lower degree \cite[Proposition~9]{LST19}. In both these cases this is essentially the formal computation $W_0 = e^{-t\Delta}W = \sum_{n \geq 0} \frac{(-t)^n}{n!}\Delta^kW$ which we see is well-defined for at least $W$ the product of a polynomial and an exponential. It would appear at first glance that exact inversion (or a meaningful description thereof) of a general $W$ is either difficult, or outright impossible, due to the smoothing properties of the heat propagator. On the latter point, the perturbation type arguments of \cite{LST19} suggest exact inversion is only required at top order, provided lower order terms can be appropriately controlled. In some sense this is not surprising; were it the case $r \in L^\infty(\bbR)$ then one expects for a typical $W$ the asymptotics are unaffected by adding or subtracting constants. 

In this paper, we shall be concerned with obtaining Weyl type results for the Riesz means of functional difference operators with potentials of the form $W(x) = \abs{x}^pe^{\abs{x}^\beta}$ for $p \geq 0$ and $\beta \in [0, 2]$ not both zero, generalising the study of polynomial and exponential potentials initiated in \cite{LST15, LST19}. The specific form of $W$ above can be replaced without much difficulty with potentials of the form $W(x) = \sum_{i=1}^n p_i(x)e^{q_i(x)}$ for $p_i, q_i$ polynomial in either $x$ or $(1 + \abs{x}^2)^{1/2}$ with $q_i$ growing at most quadratically, and in general the method presented in this paper can be adapted to the study of other spectral problems including for the classical Schr\"odinger operators $H = -\Delta + V$. 

We shall present a different proof idea based on varying the parameter $a \rightarrow \infty$ defining the gaussian $g_a$ and thus the coherent state transform. There are at least two advantages of this proof. The first is that it simplifies the analysis of the Riesz mean by simplifying the ``error" term, and the second is that it sidesteps the delicate issue of requiring the existence of a solution to the convolution equation ${W_0 \ast g_a^2} = W$ even at top order. The (leading term in the) asymptotics of the eigenvalue counting function $N(\lambda) = \#\{j \mid \lambda_j \leq \lambda\}$ follows as a corollary, thanks to a Karamata-Tauberian type theorem (see \cite[Corollary~2.3]{LST19} and \cite[Theorem~10.3]{Sim79}) relating the Riesz mean to the eigenvalue counting function. For convenience, we simplify the scaling by replacing $H_0$ with Fourier multiplication by $\cosh(\xi)$ in the sequel.

\begin{theorem}
Let $H = H_0 + W$. If $W(x) = \abs{x}^p$ for $p > 0$ then the Riesz mean of $H$ has asymptotics 
\[ \sum_{j \geq 1} (\lambda - \lambda_j)_+ \asymp_{\lambda \to \infty} \frac{4p}{p+1} \lambda^{1+1/p}\log\lambda \]
and if $W(x) = \abs{x}^pe^{\abs{x}^\beta}$ for $p \geq 0$ and $\beta \in (0, 2]$ then the Riesz mean of $H$ has asymptotics
\[ \sum_{j \geq 1} (\lambda - \lambda_j)_+ \asymp_{\lambda \to \infty} 4\lambda(\log\lambda)^{1+1/\beta}. \]
\end{theorem}

\begin{corollary}
Let $H = H_0 + W$. If $W(x) = \abs{x}^p$ for $p > 0$ then the eigenvalue counting function has asymptotics
\[ N(\lambda) \asymp_{\lambda \to \infty} 4\lambda^{1/p}\log\lambda \]
and if $W(x) = \abs{x}^pe^{\abs{x}^\beta}$ for $p \geq 0$ and $\beta \in (0, 2]$ then the eigenvalue counting function has asymptotics
\[ N(\lambda) \asymp_{\lambda \to \infty} 4(\log\lambda)^{1+1/\beta}. \]
\end{corollary}

%% file: main-results.tex
\section{Proof of main results}
We first recall the bounds on the Riesz means obtained in \cite{LST15}. Let $(\psi_j)_{j \geq 1}$ be the sequence of $L^2(\bbR)$-normalised eigenfunctions of $H$ associated to $(\lambda_j)_{j \geq 1}$. The Riesz mean is given by 
\[ \sum_{j \geq 1} (\lambda - \lambda_j)_+ = \sum_{j \geq 1} (\lambda - (H_0\psi_j, \psi_j) - (W\psi_j, \psi_j))_+. \]
It can be shown an analogue of \eqref{lst19formula} for multiplication by a function $(k, y) \mapsto T(k)$ depending only on $k$ in coherent state space holds, namely 
\[ \int_{\bbR^2} T(k)\smallabs{\tilde{\psi}(k, y)}^2dkdy = \int_\bbR (T \ast \hat{g}_a^2)(\xi)\smallabs{\hat{\psi}(\xi)}^2d\xi \]
so that, since $H_0$ is Fourier multiplication by $\cosh(\xi)$, it suffices to solve the convolution equation ${T \ast \hat{g}_a^2} = \cosh(\xi)$ to find $T(k) = e^{-a/(16\pi^2)}\cosh(k)$. Our first observation is that if $W_0$ satisfies $W \geq {W_0 \ast g_a^2}$ then following the arguments in \cite[\S2.2]{LST15} we obtain the upper bound
\begin{equation}\label{upperbound}
\begin{aligned}
    \sum_{j \geq 1} (\lambda - \lambda_j)_+ &\leq \sum_{j \geq 1} (\lambda - (T(k)\tilde{\psi}_j, \tilde{\psi}_j) - (W_0(y)\tilde{\psi}_j, \tilde{\psi}_j))_+ \vphantom{\int_{\bbR^2} (\lambda - e^{-a/(16\pi^2)}\cosh(k) - W_0(y))_+dkdy}  \\
    &= \int_{\bbR^2} (\lambda - e^{-a/(16\pi^2)}\cosh(k) - W_0(y))_+dkdy \vphantom{\int_{\bbR^2} (\lambda - e^{-a/(16\pi^2)}\cosh(k) - W_0(y))_+dkdy}.
\end{aligned}
\end{equation}
That is an exact inverse $W_0$ satisfying ${W_0 \ast g_a^2} = W$ is not necessary, and it suffices to find a suitable one-sided inverse $W_0$. Similarly, if $W_1$ satisfies ${W \ast g_a^2} \leq W_1$, then following the arguments in \cite[\S2.3]{LST15} we obtain the lower bound
\begin{equation}\label{lowerbound}
\begin{aligned}
\sum_{j \geq 1} (\lambda - \lambda_j)_+ &\geq \int_{\bbR^2} (\lambda - (\cosh \ast \, \hat{g}_a^2)(k) - (W \ast g_a^2)(y))_+dkdy \\ 
&= \int_{\bbR^2} (\lambda - e^{a/(16\pi^2)}\cosh(k) - W_1(y))_+dkdy. 
\end{aligned}
\end{equation}

\begin{remark}
The (one-sided) inverse $W_0$ plays no role in the lower bound. 
\end{remark}

We now give a loose sketch of the proof. We already know from \cite{LST15} one way to work with the coherent state transform is to exactly invert $W$, and we know from \cite{LST19} a generalisation of this is that it suffices to invert at top order, that is modulo lower order terms. In this paper, we present a proof that apparently sidesteps the inversion all together. Since convolution with $g_a^2$ is forward heat propagation for time $t(a) = 1/(4a)$, we expect as $a \rightarrow \infty$ and $t = t(a) \rightarrow 0$ that $g_a^2 \rightarrow \delta$ converges to the delta at zero in the sense of distributions. It follows by continuity of $W$ that ${W \ast g_a^2} \rightarrow W$ formally (and, rigorously, uniformly on compact subsets). On the other hand, convolution with $\hat{g}_a^2$ is forward heat propagation for time $s = s(a) = a/(16\pi^2)$ and so the coherent state multiplier $T(k)$ corresponding to $H_0$ is $e^{-s(a)}\cosh(k)$. This hides in some sense an uncertainty principle of sorts since as $a \rightarrow \infty$ we have ${W \ast g_a^2} \rightarrow W$ while $e^{-s} \rightarrow 0$, that is the error of ${W \ast g_a^2}$ with respect to $W$ can be made small provided we pay with the error of $T(k)$ with respect to $\cosh(k)$. This uncertainty principle appears in the lower bound also, since $e^s \rightarrow \infty$. 

We will be interested in finding explicit error estimates for ${W \ast g_a^2} - W$. For instance, if $W(x) = \cosh(x)$ then ${(W \ast g_a^2)(y)} = e^{1/(4a)}\cosh(y)$ so for $a$ large we expect ${W \ast g_a^2} - W \leq \sigma(a)W$ for $\sigma = \sigma(a) = o(1)$ with respect to the $a \rightarrow \infty$ limit. So let us first suppose there exists an affine estimate of the form 
\begin{equation}\label{hypothesis}
W \ast g_a^2 \leq (1 + \sigma(a))W + \tau(a) 
\end{equation}
where $\tau = \tau(a) = o(1)$. Then by \eqref{hypothesis} we have
$W_0 \ast g_a^2 \leq W$ for $W_0 \vcentcolon= (1 + \sigma)^{-1}(W - \tau)$ and hence by \eqref{upperbound} it holds
\begin{equation}\label{upperbound1}
\sum_{j \geq 1} (\lambda - \lambda_j)_+ \leq \int_{\bbR^2} (\lambda - e^{-s}\cosh(k) - (1 + \sigma)^{-1}(W(y) - \tau))_+dkdy. 
\end{equation}
Similarly, ${W \ast g_a^2} \leq (1 + \sigma)W + \tau =\vcentcolon W_1$ and hence by \eqref{lowerbound} it holds
\begin{equation}\label{lowerbound1}
\sum_{j \geq 1} (\lambda - \lambda_j)_+ \geq \int_{\bbR^2} (\lambda - e^s\cosh(k) - ((1 + \sigma)W + \tau))_+dkdy.
\end{equation}
As mentioned in the introduction, our proof simplifies the error terms in that, taking for instance as example the case $W(x) = x^{2n}$, rather than $R(y)$ a polynomial of lower degree as the error ${W \ast g_a^2} - W$, we have a linearised error in terms of $W$ which reduces analysis of both upper and lower bounds to fundamentally the same integral. There is a cost, however, in that we must define a suitable scale $a = a(\lambda)$ such that $a \rightarrow \infty$ as $\lambda \rightarrow \infty$, and we must also track estimates on $s, \sigma, \tau$ which are defined in terms of $a$. 

The choice of an affine estimate is deliberate; it includes the class of exponentials (which contribute constant prefactors) and polynomials (which contribute lower order terms). For exponential $W$ and in particular hyperbolic $W$ the estimate \eqref{hypothesis} for \emph{fixed} $a$ and therefore fixed $\sigma, \tau$ together with \eqref{upperbound1}, \eqref{lowerbound1} already implies the expected asymptotics, as the leading term of the integral $\int_{\bbR^2} (\lambda - e^{\pm s}\cosh(k) - W(y))_+dkdy$ is independent of fixed $s$ and preserved by a fixed affine transformation of $W$.

For polynomial $W$ however, the prefactor $(1 + \sigma)^{\pm 1}$ appears in the leading term and therefore some controlled limit $a \rightarrow \infty$ is required. In the applications to follow, we may be rather generous at times with the error estimates, in the sense we may simplify the analysis by using less precise estimates without any loss in the asymptotics. This approach is then tractable in that one need not track the best possible estimates (which can be rather advantageous in situations where such estimates are difficult or too laborious to obtain). 

Let us also provide a brief discussion of the asymptotics of 
\[ \int_{\bbR^2} (\lambda - C\cosh(k) - W(y))_+dkdy \]
for $C > 0$ and $W$ a differentiable real-valued function growing symmetrically and monotonically to infinity. By symmetry, the integral over $\bbR^2$ is four times the same integral taken over the first quadrant, and, by monotonicity, we may substitute either $u = \frac{C}{2}e^k$ since $\cosh(k) \geq \frac{1}{2}e^k$, or $u = Ce^k$, since $\cosh(k) \leq e^k$, depending on whether we are estimating the upper or lower bounds respectively. This reduces the analysis to the asymptotics of
\[ \int_{\bbR^2} (\lambda - Ce^k - W(y))_+dkdy \]
for $C > 0$. Substituting $v = W(y)$ on $\{y \geq 0\}$ we obtain 
\begin{equation}\label{asympint}
\begin{aligned}
\int_{\bbR^2} (\lambda - Ce^k - W(y))_+ &= 4\int_{W(0)}^{\lambda - C}\int_C^{\lambda - v} \frac{\lambda - u - v}{uW'(W^{-1}(v))}dudv \\
&= 4\int_C^{\lambda-W(0)}\int_{W(0)}^{\lambda-u} \frac{\lambda - u - v}{uW'(W^{-1}(v))}dvdu.
\end{aligned}
\end{equation}
For some choices of $W$, the leading term of this integral can be explicitly computed. For instance, if $W(y) = \cosh(y)$ or $W(y) = \abs{y}^p$ for $p > 0$ then we obtain the leading terms $\lambda\log^2\lambda$ and $\frac{p}{p+1}\lambda^{1+1/p}\log\lambda$ respectively, the latter of which we prove shortly (and the former of which was proven in \cite{LST15}). We may of course work with more general $W$. For instance, we may allow $W$ to behave arbitrarily bad in a compact neighbourhood of the origin by decomposing the region of integration into a suitable disc and its complement; we may allow for absence of symmetry by dealing with each quadrant individually (in which case one quadrant alone may determine the leading term), and we may also allow for absence of differentiability or monotonicity by, for instance, controlling $W$ above and below by a more regular function exhibiting the same growth. 

\subsection{Power potentials of form $W(x) = \abs{x}^p$ for $p > 0$}
Let $W(x) = \abs{x}^p$. To estimate ${W \ast g_a^2}$, we first estimate $W(x - y) - W(x)$ since 
\[ (W \ast g_a^2)(x) - W(x) = \int_\bbR (W(x - y) - W(x))g_a^2(y)dy. \]
If $p \geq 2$ then write $p/2 = \floor{p/2} + (p/2 - \floor{p/2}) =\vcentcolon q + r$ for $q \geq 1$ and $r \in [0, 1)$. By the binomial theorem and subadditivity of $x \mapsto \abs{x}^r$ we have
\[ \abs{x - y}^p = ((x - y)^2)^{p/2} \leq (\abs{x}^2 + 2\abs{x}\abs{y} + \abs{y}^2)^{q+r} \leq \abs{x}^p + \sum_{a, b} C_{a, b}\abs{x}^a\abs{y}^b \]
where $C_{a, b} > 0$ are constants depending on binomial coefficients, $q$, and $r$, and where $0 \leq a \leq p-1$ and $b > 0$. If $p < 2$ then 
\[ \abs{x - y}^p \leq \abs{x}^p + 2^{p/2}\abs{x}^{p/2}\abs{y}^{p/2} + \abs{y}^p. \] 
In fact if $p \leq 1$ we have the better estimate $\abs{x - y}^p \leq \abs{x}^p + \abs{y}^p$ but it suffices to proceed in the sequel with estimates of the form 
\begin{equation}\label{xpest}
\abs{x - y}^p \leq \abs{x}^p + \sum_{a, b} C_{a, b}\abs{x}^a\abs{y}^b \leq \abs{x}^p + C(1 + \abs{x}^p)\sum_\gamma \abs{y}^{\gamma}, \quad \gamma \in \bbR_{>0}
\end{equation}
for some constant $C > 0$, in which case
\begin{equation}\label{polyest}
\abs{(W \ast g_a^2)(x) - W(x)} \leq C(1 + \abs{x}^p) \abs{\int \sum_\gamma \abs{y}^{\gamma} g_a^2(y)dy}.
\end{equation}
By the triangle inequality, the postfactor on the right hand side is a sum of central absolute moments of the gaussian $g_a$ whose asymptotics are given by 
\[ \int \abs{y}^\gamma g_a^2(y)dy \leq \frac{C_\gamma}{a^{-\gamma/2}} \]
for some $C_\gamma > 0$ and $a$ sufficiently large. Thus the postfactor decays polynomially fast in $a$, in particular like $a^{-\gamma_0/2}$ where $\gamma_0 = \min\gamma$ is the smallest $\gamma$ appearing in \eqref{xpest}, and we can write 
\[ \abs{(W \ast g_a^2)(x) - W(x)} \leq C'a^{-\gamma_0/2}(1 + \abs{x}^p) \]
for some $C' > 0$. This is now of the form of the estimate \eqref{hypothesis}, i.e.
\[ W \ast g_a^2 \leq (1 + \sigma)W + \tau \]
with $\sigma = \tau = C'a^{-\gamma_0/2}$. Replacing $\tau(a)$ with $1$ for $\lambda$ and thus $a$ sufficiently large, an upper bound on the Riesz mean by \eqref{upperbound1} is 
\begin{align*}
    \sum_{j \geq 1} (\lambda - \lambda_j)_+ &\leq \int_{\bbR^2} (\lambda - e^{-s}\cosh(k) - (1 + \sigma)^{-1}(W(y) - 1))_+dkdy \\
    &= (1 + \sigma)^{-1}\int_{\bbR^2} (\lambda_1 - (1 + \sigma)e^{-s}\cosh(k) - W(y))_+dkdy \vphantom{\sum_{j \geq 1} (\lambda - \lambda_j)_+ \leq \int_{\bbR^2} (\lambda - e^{-s}\cosh(k) - (1 + \sigma)^{-1}(W(y) - 1))_+dkdy} \\
    &\leq (1 + \sigma)^{-1}\int_{\bbR^2} (\lambda_1 - C_1e^k - W(y))_+dkdy \vphantom{\sum_{j \geq 1} (\lambda - \lambda_j)_+ \leq \int_{\bbR^2} (\lambda - e^{-s}\cosh(k) - (1 + \sigma)^{-1}(W(y) - 1))_+dkdy}
\end{align*}
for $\lambda_1 = (1 + \sigma)\lambda + 1$ and $C_1 = \frac{1}{2}(1 + \sigma)e^{-s}$. Similarly, a lower bound by \eqref{lowerbound1} is
\begin{align*}
    \sum_{j \geq 1} (\lambda - \lambda_j)_+ &\geq \int_{\bbR^2} (\lambda - e^{s}\cosh(k) - ((1 + \sigma)W + \tau))_+dkdy \\
    &= (1 + \sigma)\int_{\bbR^2} (\lambda_2 - (1 + \sigma)^{-1}e^s\cosh(k) - W(y))_+dkdy \vphantom{\sum_{j \geq 1} (\lambda - \lambda_j)_+ \geq \int_{\bbR^2} (\lambda - e^{s}\cosh(k) - ((1 + \sigma)W + \tau))_+dkdy} \\
    &\geq (1 + \sigma)\int_{\bbR^2} (\lambda_2 - C_2e^k - W(y))_+dkdy \vphantom{\sum_{j \geq 1} (\lambda - \lambda_j)_+ \geq \int_{\bbR^2} (\lambda - e^{s}\cosh(k) - ((1 + \sigma)W + \tau))_+dkdy}
\end{align*}
for $\lambda_2 = (1 + \sigma)^{-1}(\lambda - 1)$ and $C_2 = (1 + \sigma)^{-1}e^s$. 

The problem is now reduced to the study of the phase space integral $\int_{\bbR^2} (\lambda - Ce^k - W(y))_+dkdy$. The computations for $W(y) = y^{2n}$ appearing in \cite[\S2.2]{LST19} can be generalised to the case $W(y) = y^p$ for $p > 0$, but for completeness we repeat the computations since we need to track the constant $C$ which depends on $s$, thus $a$, and thus $\lambda$. Note since $\sigma = o(1)$ it turns out $(1 + \sigma)^{\pm 1} = 1 + o(1)$ and therefore $\sigma$ does not appear in the leading term (so essentially the estimate on $C$ is the only meaningful one). By our comments earlier on the integral \eqref{asympint}, 
\[ \int_{\bbR^2} (\lambda - Ce^k - W(y))_+dkdy = 4\int_C^{\lambda-W(0)}\int_{W(0)}^{\lambda-u} \frac{\lambda - u - v}{uW'(W^{-1}(v))}dvdu. \]
Proceeding with $W(0) = 0$, the part of the integral with numerator $-u$ is of lower order, namely 
\begin{align*}
    \abs{\int_C^{\lambda-W(0)}\int_{W(0)}^{\lambda-u} \frac{-u}{uW'(W^{-1}(v))}dvdu} &= \int_C^\lambda\int_0^{\lambda-u} \frac{d}{dv}(W^{-1}(v))dvdu \\
    &= \int_C^\lambda W^{-1}(\lambda - u)du \vphantom{\abs{\int_C^{\lambda-W(0)}\int_{W(0)}^{\lambda-u} \frac{-u}{uW'(W^{-1}(v))}dvdu}} 
\end{align*}
which is bounded above by $\lambda W^{-1}(\lambda) = \lambda^{1+1/p}$, provided $C < \lambda$. The leading term appears in the part of the integral with numerator $\lambda - v$. As done in \cite{LST15, LST19}, we simplify the integral by rescaling $u' = u/\lambda$ and $v' = v/\lambda$ to reduce to 
\begin{align*} 
    \int_C^{\lambda-W(0)}\int_{W(0)}^{\lambda-u} \frac{\lambda-v}{uW'(W^{-1}(v))}dvdu &= \frac{1}{p} \int_C^\lambda\int_0^{\lambda-u} \frac{\lambda-v}{uv^{(p-1)/p}}dvdu \\
    &= \frac{\lambda^{1+1/p}}{p}\int_{C/\lambda}^1 \int_0^{1-u} \frac{1-v}{uv^{(p-1)/p}}dvdu \\
    &= \frac{\lambda^{1+1/p}}{p+1}\int_{C/\lambda}^1 \frac{(1-u)^{1/p}(p+u)}{u}du. 
\end{align*}
Integrating by parts, 
\begin{align*}
    \int_{C/\lambda}^1 \frac{(1-u)^{1/p}(p+u)}{u}du &= (\log \lambda/C)(1 - C/\lambda)^{1/p}(p + C/\lambda) \vphantom{\frac{p+1}{p}\int_{C/\lambda}^1 u(1-u)^{1/p-1}\log(u)du} \\
    &+ \frac{p+1}{p}\int_{C/\lambda}^1 u(1-u)^{1/p-1}\log(u)du. \vphantom{\int_{C/\lambda}^1 \frac{(1-u)^{1/p}(p+u)}{u}du}
\end{align*}
The first addend contributes the leading term while the other addend is $\cO(1)$ independently of $C$ since $\int_0^1 u(1 - u)^{1/p-1}\log(u)du < \infty$. 

To summarise the analysis, provided $C = C_1 < \lambda_1 = \lambda$, an upper bound on the Riesz mean is given by 
\begin{equation}\label{ubasymppoly}
4(1 + \sigma)^{-1}\frac{\lambda_1^{1+1/p}}{p+1}(\log \lambda_1/C_1)(1 - C_1/\lambda_1)^{1/p}(p + C_1/\lambda_1) + \cO(\lambda_1^{1+1/p}) + \cO(1)
\end{equation}
and hence it suffices to choose the scale $a(\lambda) = \log\log\lambda$ from which it follows $\sigma = o(1)$, $\lambda_1 = (1 + o(1))\lambda + \cO(1)$, and $C = \cO((\log\lambda)^{-1/(16\pi^2)}) = o(\lambda)$. Expanding \eqref{ubasymppoly} according to this choice of scale gives the expected asymptotics according to the theorem. Computations for the lower bound are essentially identical, with $(1 + \sigma)$, $\lambda_2 = (1 + o(1))\lambda + \cO(1)$, and $C_2=\cO((\log\lambda)^{1/(16\pi^2)}$ replacing $(1 + \sigma)^{-1}$, $\lambda_1$, and $C_1$ respectively. 

\subsection{Exponential potentials of form $W(x) = \abs{x}^pe^{\abs{x}^\beta}$ for $p \geq 0$ and $\beta \in (0, 2]$}
Let $W(x) = \abs{x}^pe^{\abs{x}^\beta}$. It turns out in the subexponential and exponential cases $\beta \leq 1$ the idea of the previous section is sufficient, in that it can be used without significant modification to obtain the expected asymptotics, thanks to the inequality $\abs{x - y}^\beta \leq \abs{x}^\beta + \abs{y}^\beta$. To handle superexponential growth however, we play the same game but with a different version of the estimate \eqref{hypothesis} that can handle the fact the heat propagator acting on (the reciprocal of) a gaussian dilates its variance. To this end, first note 
\[ \abs{x - y}^\beta \leq \abs{x}^\beta + 2^{\beta/2}\abs{y}^{\beta/2}\abs{y}^{\beta/2} + \abs{y}^\beta \leq (1 + \epsilon)\abs{x}^\beta + (1 + K(\epsilon))\abs{y}^\beta \]
by Young's inequality and for a scale $\epsilon = \epsilon(\lambda) \rightarrow 0$ to be determined later, and $K = K(\epsilon) = \cO(\epsilon^{-1})$. We introduce a dilation on the argument of $W$ by $\mu = \mu(\epsilon) = (1 + \epsilon)^{1/\beta}$. Denoting the summation in \eqref{polyest} with $\Sigma = \Sigma(y) \vcentcolon= \sum_\gamma \abs{y}^\gamma$, 
\begin{align*}
    W(x - y) - \mu^{-p}W(\mu x) &\leq (\abs{x}^p(1 + C\Sigma) + C\Sigma)e^{(1+\epsilon)\abs{x}^\beta}e^{(1+K)\abs{y}^\beta} - \abs{x}^pe^{(1+\epsilon)\abs{x}^\beta} \\
    &\leq \abs{x}^pe^{(1+\epsilon)\abs{x}^\beta}((1 + C\Sigma)e^{(1+K)\abs{y}^\beta}-1) + C\Sigma e^{(1+\epsilon)\abs{x}^\beta}e^{(1+K)\abs{y}^\beta} \\
    &\leq \abs{x}^pe^{(1+\epsilon)\abs{x}^\beta}((1 + 2C\Sigma)e^{(1+K)\abs{y}^\beta}-1) + e^2C\Sigma e^{(1+K)\abs{y}^\beta} \\ 
    &= \mu^{-p}W(\mu x)((1 + 2C\Sigma)e^{(1+K)\abs{y}^\beta}-1) + e^2C\Sigma e^{(1+K)\abs{y}^\beta}
\end{align*}
where in the third inequality we used $e^{(1 + \epsilon)\abs{x}^\beta} \leq \abs{x}^pe^{(1 + \epsilon)\abs{x}^\beta} + e^2$, valid for $\epsilon \leq 1$ and thus in the $\lambda \rightarrow \infty$ limit. In the setting of $\beta > 0$ it turns out constant prefactors, i.e. $(1 + \sigma)^\pm$, do not appear in the leading term and, like in the polynomial setting of the previous chapter, constant additional terms, i.e. $\pm\tau$, also do not appear. Thus we can be generous with our estimates to simplify the proof. In particular, choosing scales for $a$ and $\epsilon$ such that
\begin{equation}\label{Kscale}
    \int ((1 + 2C\Sigma)e^{(1 + K)\abs{y}^\beta}-1)g_a^2(y)dy \rightarrow 0 \Mand \int e^2C\Sigma e^{(1+K)\abs{y}^\beta}g_a^2(y)dy \rightarrow 0 
\end{equation}
as $\lambda \rightarrow \infty$, for $\lambda$ and thus $a$ sufficiently large it holds ${(W \ast g_a^2)(x)} - \mu^{-p}W(\mu x) \leq (2 - \mu^{-p})W(\mu x) + 2$ and hence 
\begin{equation}\label{secondhypothesis}
W \ast g_a^2 \leq 2(W(\mu x) + 1). 
\end{equation}
We will need to track the estimates on $\epsilon$ since the scale $\mu$ dilating the argument of $W$ does affect the leading term and we need to argue we can choose a scale for $\epsilon$ such that $\mu \rightarrow 1$ preserves the expected asymptotics and ensures \eqref{Kscale} holds; if $\epsilon \rightarrow 0$ very rapidly then $K \rightarrow \infty$ and $\beta = 2$ may violate \eqref{Kscale}. We also need an estimate in the opposite direction since \eqref{secondhypothesis} only provides the lower bound. If $\tilde{W}(x) = W(cx)$ for $c > 0$ it is readily checked ${(\tilde{W} \ast g_a^2)(x)} = {(W \ast g_{a/c^2})(cx)} \leq 2(W(\mu c x) + 1)$ for $\lambda$ and thus $a/c^2$ sufficiently large. Since $\mu \rightarrow 1$ it follows $W_0(x) = \frac{1}{2}W(\mu^{-1}x) - 1$ will provide the upper bound. 

Let us complete the proof by showing we can choose appropriate scales for $a$ and $\epsilon$ ensuring the previous statements hold. An upper bound on the Riesz mean is given by 
\begin{align*}
    \sum_{j \geq 1} (\lambda - \lambda_j)_+ &\leq \int_{\bbR^2} \left(\lambda - e^{-s}\cosh(k) - \left(\frac{1}{2}W(\mu^{-1}y) - 1\right)\right)_+dkdy \\
    &= \frac{1}{2}\int_{\bbR^2} (\lambda_1 - 2e^{-s}\cosh(k) - W(\mu^{-1}y))_+dkdy \vphantom{\int_{\bbR^2} \left(\lambda - e^{-s}\cosh(k) - \left(\frac{1}{2}W(\mu^{-1}y) - 1\right)\right)_+dkdy} \\
    &\leq \frac{\mu}{2}\int_{\bbR^2} (\lambda_1 - C_1e^k - W(y))_+dkdy \vphantom{\int_{\bbR^2} \left(\lambda - e^{-s}\cosh(k) - \left(\frac{1}{2}W(\mu^{-1}y) - 1\right)\right)_+dkdy}
\end{align*}
for $\lambda_1 = 2(\lambda + 1)$ and $C_1 = e^{-s}$. Similarly, a lower bound is given by
\[ \sum_{j \geq 1} (\lambda - \lambda_j)_+ \geq \frac{2}{\mu} \int_{\bbR^2} (\lambda_2 - C_2e^k - W(y))_+dkdy \]
for $\lambda_2 = \frac{\lambda}{2} - 1$ and $C_2 = \frac{1}{2}e^s$. The problem is again now reduced to the study of the phase space integral $\int_{\bbR^2} (\lambda - Ce^k - W(y))_+dkdy$ with $\lambda = \lambda_i$ and $C = C_i$ for $i \in \{1, 2\}$. Note since $\mu = 1 + o(1) = \mu^{-1}$ for any choice of scale $\epsilon \rightarrow 0$, it makes no appearance in the leading term so that, again, essentially the estimate on $C$, depending on $s$ through $a$, is the only meaningful one, and hence the only question is whether we can choose $a \rightarrow \infty$ and $\epsilon \rightarrow 0$ correctly. As before,
\[ \int_{\bbR^2} (\lambda - Ce^k - W(y))_+dkdy = 4\int_C^{\lambda-W(0)}\int_{W(0)}^{\lambda-u} \frac{\lambda - u - v}{uW'(W^{-1}(v))}dvdu \] 
and again the part with numerator $-u$ 
\[ \abs{\int_C^\lambda\int_0^{\lambda-u} \frac{-u}{uW'(W^{-1}(v))}} = \cO(\lambda W^{-1}(\lambda)) = \cO(\lambda (\log \lambda)^{1/\beta}) \]
is of lower order. For the part with numerator $\lambda - v$, there is a closed form expression for the inverse of $W$ in terms of the Lambert-$W$ function. This integral is still rather unwieldy however, but we can simplify by cutting away from a neighbourhood of the origin $\{0 \leq v \leq L\}$ for $L \gg 1$ sufficiently large independently of $\lambda$ and insert explicit asymptotics for $W^{-1}$. The integral over $\{0 \leq v \leq L\}$ is of lower order because 
\begin{align*}
    \int_C^{\lambda - L}\int_0^L \frac{\lambda - v}{uW'(W^{-1}(v))}dvdu &\leq \lambda \int_C^{\lambda - L} \int_0^L \frac{1}{u(W'(W^{-1}(v))}dvdu  \\
    &= \lambda\int_C^{\lambda - L} \frac{du}{u} \int_0^L \frac{d}{dv}W^{-1}(v)dv \\ 
    &= \cO(\lambda \log \lambda/C),\vphantom{\lambda\int_C^{\lambda - L} \frac{du}{u} \int_0^L \frac{d}{dv}W^{-1}(v)dv}
\end{align*}
provided $C$ is appropriately chosen. As for the integral over $\{L \leq v \leq \lambda - u\}$, integrating by parts,
\[ \int_C^{\lambda-L}\int_L^{\lambda-u} \frac{\lambda-v}{uW'(W^{-1}(v))}dvdu = \int_C^{\lambda-L} \left(\frac{(\lambda-v)W^{-1}(v)}{u}\Big\rvert_{v=L}^{\lambda-u} + \int_L^{\lambda-u} \frac{W^{-1}(v)}{u}dv\right)du, \]
the first addend is 
\[ \int_C^{\lambda-L} \left(W^{-1}(\lambda-u) + \frac{\lambda-L}{u}W^{-1}(L)\right)du = \cO(\lambda (\log\lambda)^{1/\beta}) + \cO(\lambda\log\lambda/C) \]
which again is of lower order. In the second addend we may insert the large $v$ asymptotic 
\[ W^{-1}(v) = (\log v)^{1/\beta} + o((\log v)^{1/\beta}) \]
which reduces the problem to the study of 
\[ \int_C^{\lambda-L}\int_L^{\lambda-u} \frac{(\log v)^{1/\beta}}{u}dvdu. \]

\begin{remark}
    An explicit estimate on $o((\log v)^{1/\beta})$ is $\cO((\log v)^{1/\beta-\delta})$ for some $\delta > 0$, in which case a crude estimate on the resultant integral is $\cO(\lambda (\log \lambda)^{1/\beta-\delta}\log(\lambda/C))$.
\end{remark}

For integral $1/\beta \in \bbZ_{\geq 1}$ the integral can be expressed in closed form in terms of logarithms and polylogarithms from which the leading term can be extracted directly. To avoid use of such machinery and to generalise for nonintegral $1/\beta \notin \bbZ_{\geq 1}$, we may proceed as follows. Integrating by parts with $r = 1/\beta$, 
\[ \int_C^{\lambda-L} \int_L^{\lambda-u} \frac{(\log v)^r}{u}dvdu = \int_C^{\lambda-L}\left(\frac{v(\log v)^r}{u}\Big\rvert_{v=L}^{\lambda-u} - r\int_L^{\lambda-u} \frac{(\log v)^{r-1}}{u}dv\right)du \]
and the second addend is either $\cO(\lambda(\log \lambda)^{r-1}\log(\lambda/C))$ or $\cO(\lambda\log \lambda/C)$ depending on whether $r \geq 1$ or $r < 1$ respectively; in both cases it is of lower order. For the first addend, evaluation at the lower terminal yields $\cO(u^{-1})$ and $\int_C^{\lambda-L} \cO(u^{-1})du = \cO(\log \lambda/C)$ is of lower order, whereas evaluation at the upper terminal contributes 
\[ \int_C^{\lambda-L} \frac{(\lambda - u)(\log(\lambda - u))^r}{u}du = \lambda \int_{C/\lambda}^{1-L/\lambda} \frac{(1-u)(\log \lambda + \log(1-u))^r}{u}du \]
after rescaling. Irrespective of whether the asymptotics of this integral can be computed exactly, as far as upper and lower bounds are concerned, since $r > 0$ we may insert the first part of the estimate \eqref{xpest}; the integral corresponding to the top order term $(\log \lambda)^r$ contributes the leading term since we may compute the integral explicitly
\begin{align*}
    \lambda(\log \lambda)^r\int_{C/\lambda}^{1-L/\lambda} \frac{1-u}{u}du &= \lambda(\log\lambda)^r \cdot -u + \log u\Big\rvert_{u=C/\lambda}^{1-L/\lambda} \\ 
    &= \lambda(\log \lambda)^r\log(\lambda/C) + \cO(C(\log\lambda)^r) + \cO(\lambda(\log \lambda)^r). \vphantom{\lambda(\log\lambda)^r \cdot -u + \log u\Big\rvert_{u=C/\lambda}^{1-L/\lambda}}
\end{align*}
The remaining terms are integrals of the form 
\[ \lambda (\log \lambda)^q \int_{C/\lambda}^{1-L/\lambda} \frac{(1 - u)\abs{\log(1 - u)}^{q'}du}{u} \]
for $0 \leq q < r$ and $0 < q' < r$ modulo constants and therefore are all of lower order $\cO(\lambda(\log \lambda)^q)$ since $\int_0^1 u^{-1}(1 - u)\abs{\log(1 - u)}^{q'}du < \infty$ for any $q' > 0$. Putting this all together, provided $C = C_1 < \lambda_1 = \lambda$, an upper bound on the Riesz mean is given by
\begin{equation}\label{ubasympexp}
    \frac{\mu}{2}\lambda_1(\log \lambda_1)^r(1 - C_1/\lambda_1)\log \lambda_1/C_1 + \text{lower order terms.}
\end{equation}
To complete the proof it suffices to make a choice of scales $a$ and $\epsilon$. As mentioned earlier, any (vanishing) scale for $\epsilon$ gives ${\mu = 1 + o(1)}$ and so the only consideration is a scale for which the convergence \eqref{Kscale} holds. As before, we may take $a(\lambda) = \log\log\lambda$ and $\epsilon = (\log\log\log\lambda)^{-1}$ so that $K(\epsilon) = \cO(\log\log\log\lambda)$ does not blow up too fast and violate \eqref{Kscale}. It is readily checked expanding \eqref{ubasympexp} according to these choices of scales gives the expected asymptotics according to the theorem. Computations for the lower bound are again identical; the only important difference is the expansion of ${(\log \lambda + \log(1 - u))^r \leq (\log \lambda)^r + R}$, where $R$ is a remainder involving products of lower order powers of $\log \lambda$ and powers of $\abs{\log(1 - u)}$, uses an upper bound (preserving the correct leading term) in contrast to the computations of the previous section where the integral containing the leading term was computed exactly. For the lower bound we may write instead ${(\log \lambda + \log(1 - u))^r \geq (\log \lambda)^r - \tilde{R}}$ where $\tilde{R}$ is another remainder of lower order and proceed as before. 